\pdfoutput=1
\documentclass[11pt,a4paper]{article}

\usepackage[utf8]{inputenc}
\usepackage[T1]{fontenc}
\usepackage{lmodern}
\usepackage[english]{babel}
\usepackage{microtype}

\usepackage[margin=2.5cm]{geometry}

\usepackage{amsmath,amssymb,amsthm,mathtools}
\usepackage{bm}
\usepackage{xfrac}
\usepackage{mathrsfs}

\usepackage{graphicx}
\usepackage{caption}
\usepackage{tikz}
\usetikzlibrary{arrows.meta,decorations.pathreplacing}

\usepackage[round]{natbib}
\usepackage[dvipsnames]{xcolor}
\usepackage[colorlinks=true,linkcolor=MidnightBlue,citecolor=OliveGreen,
            urlcolor=MidnightBlue]{hyperref}
\usepackage[capitalise,nameinlink]{cleveref}

\theoremstyle{plain}
\newtheorem{theorem}{Theorem}[section]
\newtheorem{lemma}[theorem]{Lemma}
\newtheorem{proposition}[theorem]{Proposition}
\newtheorem{corollary}[theorem]{Corollary}

\theoremstyle{remark}
\newtheorem{remark}[theorem]{Remark}

\title{Metallic Circles: Concyclic Tangency Points of\\ Ford Circles over \texorpdfstring{$k$}{k}-Fibonacci Quotients}

\author{%
  Nicolás Piña\thanks{Faculty of
  Engineering and Architecture, Universidad Central de Chile. E-mail:
  \href{mailto:nicolas.pina@ucentral.cl}{\texttt{nicolas.pina@ucentral.cl}}}
}

\date{\today}

\begin{document}

\maketitle

\begin{abstract}
Extending a result of \citet{kocik2020fibonacci}, we show that the tangency points of the Ford circles built
over quotients of $k$-Fibonacci numbers are concyclic: they lie on a circle centred on
the horizontal axis, of radius $\sqrt{k^{2}+4}/2$, meeting it at the reciprocal of the
metallic ratio $\sigma_{k}$ and at $-\sigma_{k}$. Concyclicity turns out to be a
quadratic identity of the sequence read in coordinates. In the Poincaré half-plane, this
circle is a geodesic along which the tangency points advance with constant step
$2\ln\sigma_k$.
\end{abstract}



\section{Introduction}\label{sec:intro}

In a short and interesting paper, \citet{kocik2020fibonacci} shows that the Fibonacci sequence hides a circle. If $F_{n}$ denotes the $n$-th Fibonacci number, the points
\begin{equation*}
    \mathcal{P}_{n}=\left( \frac{F_{2n}}{F_{2n+1}},\frac{1}{F_{2n+1} } \right), \qquad n\in \mathbb{Z},
\end{equation*}
are concyclic; more precisely, they lie on the circle with centre $(-\frac{1}{2},0)$ and radius $\frac{\sqrt{5}}{2}$, which meets the $X$-axis at the points $(\varphi^{-1},0)$ and $(-\varphi,0)$, where $\varphi$ is the golden ratio. What is remarkable about this result is where these points come from: each of them is \textit{the tangency point of two Ford circles} built over fractions $\frac{F_{n}}{F_{n+1}}$. In fact, Kocik exhibits \textit{four} concyclic sequences of this kind, labelled (A)--(D), associated with four circles of common radius $\frac{\sqrt5}{2}$ and centres $-\frac{1}{2},\frac{1}{2},\frac{3}{2},-\frac{3}{2}$; they are obtained from one another by the integer translations $x\mapsto x\pm1$, which preserve the Ford configuration. The sequence $\mathcal{P}_{n}$ above is the first of them, and it is our starting point. Independently, \citet{Sporn2022} studied circles passing through infinitely many points whose coordinates are quotients of Fibonacci numbers. Motivated by these works, the present paper answers the following question:

\begin{center}
\textit{Is the golden ratio an exceptional case, or is it part of a family?}
\end{center}

The natural candidates are the metallic ratios $\sigma_k=\frac{k + \sqrt{k^{2}+4}}{2}$ \citep{spinadel1999metallic} (for $k=1$, the golden ratio; for $k=2$, the silver ratio; and so on), together with their associated $k$-Fibonacci sequences. In \cref{sec:metallic} we prove that for each $k$ there is a circle $C_{k}$ containing all the tangency points of the Ford chain associated with the $k$-Fibonacci sequence, and meeting the $X$-axis at the points $(\sigma_{k}^{-1},0)$ and $(-\sigma_{k},0)$. Kocik's result is thus the special case $k=1$ of what happens for the whole metallic family and Ford circles.

The proof reveals that concyclicity is no mere geometric coincidence: the equation of $C_{k}$ is nothing but the identity $f_{2n+1}^{2}-kf_{2n}f_{2n+1}-f_{2n}^{2}=1$ of the $k$-Fibonacci sequence $\{f_{n}\}$, read in Cartesian coordinates. Hyperbolic geometry explains why this happens. The upper half of $C_{k}$ is a geodesic of the Poincaré half-plane and is the axis of a glide reflection that sends each tangency point to the next one. Consequently, the points $\mathcal P_n$ are equally spaced along the geodesic, at hyperbolic distance $2\ln\sigma_k$ from one another, and they travel from one foot of $C_{k}$ to the other.

\section{Preliminaries}\label{sec:prelim}

\subsection{Ford circles}\label{sec:Ford}
Given a reduced fraction $p/q$ with $q>0$, the \textit{Ford circle} $K[p,q]$
\citep{ford1938fractions} is the circle with centre $\left(\frac{p}{q},
\frac{1}{2q^{2}}\right)$ and radius $\frac{1}{2q^{2}}$; in particular, it is tangent to the
$X$-axis at $x=p/q$. We also adopt the convention of allowing $q=0$: the pair $(p,q)=(1,0)$
represents the point at infinity $x=1/0$, and $K[1,0]$ denotes the line $y=1$, the
\textit{Ford circle at infinity}. This convention is consistent with
Proposition~\ref{prop:ford}: $\det\left(\begin{smallmatrix} p & 1\\ q & 0
\end{smallmatrix}\right)=-q$ equals $\pm1$ if and only if $q=1$, and indeed the line $y=1$
is tangent precisely to the circles $K[p,1]$, all of radius $\frac12$ and with
centre at height $\frac12$. The same holds for~\eqref{eq:tangencyK}, which for
$(p,q)=(1,0)$ and $(r,s)=(r,1)$ yields the point $(r,1)$.

\begin{proposition}[\citealp{ford1938fractions}]\label{prop:ford}
Let $K[p,q]$ and $K[r,s]$ be two distinct Ford circles. If
$\det\left( \begin{smallmatrix} p & r \\ q & s \end{smallmatrix} \right)=\pm1$,
then $K[p,q]$ and $K[r,s]$ are tangent; otherwise, they are disjoint.
\end{proposition}

\begin{proposition}[\citealp{kocik2020fibonacci}]\label{prop:tangency}
If $K[p,q]$ and $K[r,s]$ are tangent, their tangency point is
\begin{equation}\label{eq:tangencyK}
\left( \frac{pq + rs}{q^{2} + s^{2}},\; \frac{1}{q^{2} + s^{2}} \right).
\end{equation}
\end{proposition}

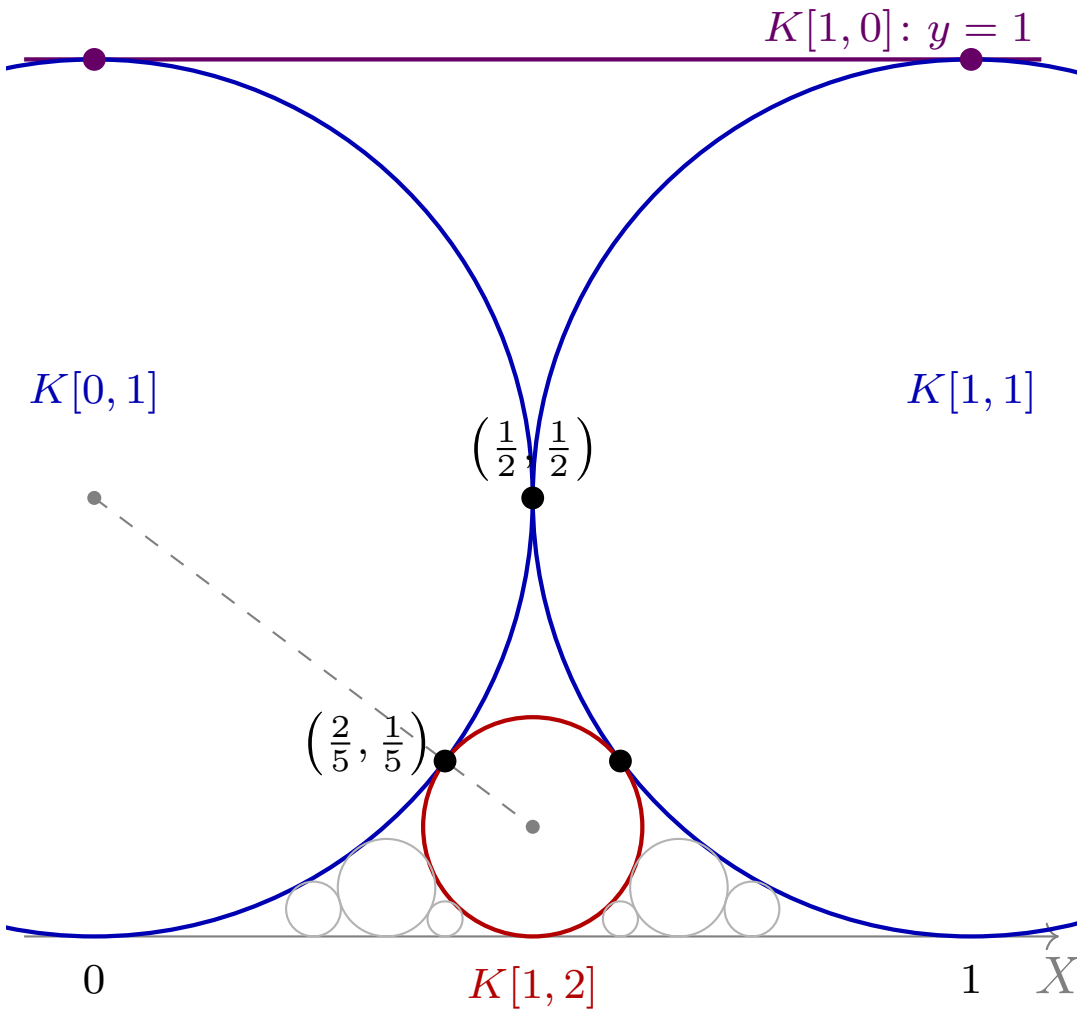
\begin{figure}[htbp]
\centering
\begin{tikzpicture}[scale=5.8]
  \clip (-0.1,-0.12) rectangle (1.12,1.1);
  \draw[->,gray] (-0.08,0) -- (1.1,0) node[below,gray] {\small $X$};
  \draw[violet!80!black,thick] (-0.08,1) -- (1.08,1) node[above left,inner sep=1pt] {\scriptsize $K[1,0]\colon y=1$};
  \draw[blue!70!black,thick] (0,0.5) circle (0.5);
  \draw[blue!70!black,thick] (1,0.5) circle (0.5);
  \draw[red!70!black,thick] (0.5,0.125) circle (0.125);
  \draw[gray!60] (0.3333,0.0556) circle (0.0556);
  \draw[gray!60] (0.6667,0.0556) circle (0.0556);
  \draw[gray!60] (0.4,0.02) circle (0.02);
  \draw[gray!60] (0.6,0.02) circle (0.02);
  \draw[gray!60] (0.25,0.03125) circle (0.03125);
  \draw[gray!60] (0.75,0.03125) circle (0.03125);
  \draw[dashed,gray] (0,0.5) -- (0.5,0.125);
  \fill[gray] (0,0.5) circle (0.008);
  \fill[gray] (0.5,0.125) circle (0.008);
  \fill[black] (0.4,0.2) circle (0.013);
  \fill[black] (0.6,0.2) circle (0.013);
  \fill[black] (0.5,0.5) circle (0.013);
  \fill[violet!80!black] (0,1) circle (0.013);
  \fill[violet!80!black] (1,1) circle (0.013);
  \node[left,inner sep=2pt] at (0.4,0.22) {\scriptsize $\left(\tfrac25,\tfrac15\right)$};
  \node[above,inner sep=3pt] at (0.5,0.5) {\scriptsize $\left(\tfrac12,\tfrac12\right)$};
  \node[blue!70!black] at (0,0.62) {\scriptsize $K[0,1]$};
  \node[blue!70!black] at (1,0.62) {\scriptsize $K[1,1]$};
  \node[red!70!black,below] at (0.5,-0.01) {\scriptsize $K[1,2]$};
  \node[below] at (0,-0.01) {\scriptsize $0$};
  \node[below] at (1,-0.01) {\scriptsize $1$};
\end{tikzpicture}
\caption{Ford circles over $\frac01$, $\frac11$ and $\frac12$ (in grey, some with
larger denominators). Since $\det\left(\begin{smallmatrix}0&1\\1&2\end{smallmatrix}\right)=-1$,
the circles $K[0,1]$ and $K[1,2]$ are tangent, and by~\eqref{eq:tangencyK} they touch at
$\left(\frac{0\cdot1+1\cdot2}{1^2+2^2},\frac{1}{1^2+2^2}\right)=\left(\frac25,\frac15\right)$,
on the segment joining their centres. The line $y=1$ plays the role of the circle $K[1,0]$.}
\label{fig:tangency}
\end{figure}

Proposition~\ref{prop:ford} is proved in \citet{ford1938fractions}; see also
\citet[\S\,5.5]{apostol1997modular}. For~\eqref{eq:tangencyK} it suffices to note that the tangency point of two externally tangent circles with centres $c_{1}$ and $c_{2}$ and radii $\rho_{1}$ and $\rho_{2}$ is $\frac{\rho_2c_1+\rho_1c_2}{\rho_1+\rho_2}$. Taking $\rho_1=\frac{1}{2q^2}$ and $\rho_2=\frac{1}{2s^2}$ and multiplying numerator and denominator by $2q^2s^2$, one obtains~\eqref{eq:tangencyK} directly. This form, symmetric in both circles, is due to Kocik; \citet[\S\,5.5]{apostol1997modular} gives an equivalent expression that is not symmetric. \Cref{fig:tangency} illustrates the formula.

\subsection{The \texorpdfstring{$k$}{k}-Fibonacci sequence and metallic ratios}\label{sec:kFib}

The Fibonacci sequence $\{F_n\}$ is defined by the recurrence relation
\[F_{n+1}=F_{n}+F_{n-1},\]
together with the initial conditions $F_{0}=0$ and $F_{1}=1$.

As one generalisation (among many) of the Fibonacci sequence, we consider, for
each integer $k\geqslant1$, the $k$-Fibonacci sequence \citep{falcon2007kfibonacci},
which we denote by $\{f_{n}\}$ and which is defined by
\begin{equation}
    f_{n+1}=kf_{n}+f_{n-1}, \qquad f_{0}=0,\quad f_{1}=1,
\end{equation}
so that $f_{n}=F_{n}$ when $k=1$. Terms with negative index are obtained by running the recurrence backwards, $f_{n-1}=f_{n+1}-kf_n$, where $n$ is a non-positive integer; by induction, $f_{-m}=(-1)^{m+1}f_m$. As for the Fibonacci sequence, the recurrence matrix $M=\left( \begin{smallmatrix} k & 1 \\ 1 & 0 \end{smallmatrix} \right)$ satisfies
\begin{equation*}
    \det M=-1, \qquad M^{n}=\begin{pmatrix}
f_{n+1} & f_{n} \\
f_{n} & f_{n-1}
\end{pmatrix}, \qquad n\in \mathbb{Z}.
\end{equation*}
Moreover, by Binet's formula,
\begin{equation}
    \lim_{n\to\infty}\dfrac{f_{n+1}}{f_{n}}=\sigma_{k},
\end{equation}
where $\sigma_{k}$ is the metallic ratio, the positive root of the quadratic equation $x^{2}-kx-1=0$. Explicitly,
\begin{equation}\label{eq:sigma_k}
\sigma_k = \frac{k + \sqrt{k^{2} + 4}}{2},
\end{equation}
whose conjugate $\sigma_{k}'$ is
\begin{equation}\label{eq:sigma_kconj}
\sigma_k' = \frac{k - \sqrt{k^{2} + 4}}{2}.
\end{equation}
Since $\sigma_k\sigma_k'=-1$, we have $\sigma_k^{-1}=-\sigma_k'$.

\subsubsection*{Special cases: golden and silver}

For the golden case we take $k=1$. The recurrence $f_{n+1}=f_{n}+f_{n-1}$ generates the Fibonacci sequence $0,1,1,2,3,5,8,\dots$ and the limit of the quotient is the golden ratio
\begin{equation*}
    \sigma_{1}=\varphi=\frac{1+\sqrt{5}}{2},\qquad \sigma_{1}'=\frac{1-\sqrt{5}}{2}=-\varphi^{-1},
\end{equation*}
which are the roots of $x^{2}-x-1=0$. The fractions $\frac{1}{1}, \frac{2}{1},\frac{3}{2}, \frac{5}{3},\frac{8}{5},\dots$ are rational approximations of $\varphi$ given by the expansion $\varphi=[1; 1, 1,\dots]$, and their reciprocals $f_{n}/f_{n+1}$ approximate $\varphi^{-1}=[0; 1, 1,\dots]$, where
\[
  [a_0; a_1, a_2, \dots] := a_0 + \cfrac{1}{a_1 + \cfrac{1}{a_2 + \cdots}}.
\]
The Ford circles $K[f_{n-1},f_{n}]$ built over the latter therefore form a chain of consecutively tangent circles whose points of contact with the $X$-axis converge to $\varphi^{-1}$.

For the silver case we take $k=2$. The recurrence $f_{n+1}=2f_{n}+f_{n-1}$ generates the Pell sequence $0,1,2,5,12,29,70,\dots$ and the limit of the quotient is the silver ratio
\begin{equation*}
\sigma_{2}=\delta_{\text{Ag}}=1+\sqrt{2},\qquad \sigma_{2}'=1-\sqrt{2}=-\delta_{\text{Ag}}^{-1},
\end{equation*}
which are the roots of $x^{2}-2x-1=0$. The fractions $\frac{2}{1}, \frac{5}{2},\frac{12}{5}, \frac{29}{12},\frac{70}{29},\dots$ are rational approximations of $\delta_{\text{Ag}}$ given by the expansion $\delta_{\text{Ag}}=[2; 2, 2,\dots]$, and their reciprocals $f_{n}/f_{n+1}$ approximate $\delta_{\text{Ag}}^{-1}=[0; 2, 2,\dots]$. The Ford circles $K[f_{n-1},f_{n}]$ built over the latter therefore form a chain of consecutively tangent circles whose points of contact with the $X$-axis converge to $\delta_{\text{Ag}}^{-1}$. Apart from a few initial circles, they differ from those of the Fibonacci chain (see \cref{fig:metallic}).

Explicitly, the fractions $f_{n-1}/f_{n}$ appearing in
\cref{fig:metallic} are: for $k=1$,
$-\tfrac32,-\tfrac21,-\tfrac11,\tfrac01,\tfrac11,\tfrac12,\tfrac23$; for $k=2$,
$-\tfrac52,-\tfrac21,\tfrac01,\tfrac12,\tfrac25$ (the terms with negative index are
obtained from $f_{-m}=(-1)^{m+1}f_{m}$). Both chains share the line $y=1$ ---which
plays the role of the circle $K[1,0]$--- and the circles $K[-2,1]$, $K[0,1]$ and $K[1,2]$.

\subsection{Five identities}\label{sec:identities}

The representation by powers of $M$ is merely a notational device: it allows us to rewrite the algebraic identities of the sequence $\{f_n\}$ as elementary matrix manipulations of $M^n$. In this subsection we gather the five identities that will be used in the subsequent results. The first three are classical for the Fibonacci sequence and extend unchanged to its $k$-Fibonacci generalisation; the fourth makes explicit the role of the parameter $k$ in the sequences under study; the fifth is d'Ocagne's identity.

\begin{lemma}\label{lem:identities}
For all $a,b,n\in\mathbb{Z}$,
\begin{align}
 f_{2n}&=f_n\left(f_{n-1}+f_{n+1}\right),  \label{eq:add1} \\
      f_{2n+1}&=f_n^2+f_{n+1}^2, \label{eq:add2}\\
f_{n-1}f_{n+1}-f_n^2&=(-1)^n, \label{eq:cassini}\\
f_{n+1}^2-k\,f_nf_{n+1}-f_n^2&=(-1)^n, \label{eq:inv}\\
f_{a}f_{b+1}-f_{a+1}f_{b}&=(-1)^{b}f_{a-b}. \label{eq:docagne}
\end{align}
\end{lemma}

\begin{proof}
All five identities follow directly from $M^{m+n}=M^mM^n$. Expanding the product on the right and comparing entries yields the addition formula
\[
  f_{m+n}=f_mf_{n+1}+f_{m-1}f_n,
\]
together with its companion $f_{m+n+1}=f_{m+1}f_{n+1}+f_mf_n$; setting $m=n$ gives~\eqref{eq:add1} and~\eqref{eq:add2}. Taking determinants of $M^{n}$ and using $\det M=-1$, we get
\[
  (-1)^{n}=\det M^n=f_{n+1}f_{n-1}-f_n^2,
\]
which is precisely Cassini's identity~\eqref{eq:cassini}. Substituting $f_{n-1}=f_{n+1}-kf_n$ into~\eqref{eq:cassini} and rearranging gives~\eqref{eq:inv}.

Finally, comparing entries in $M^{a-b}=M^{a}\left(M^{b}\right)^{-1}$ yields
d'Ocagne's identity~\eqref{eq:docagne}; see also
\citet[Ch.~III]{vajda1989fibonacci}.
\end{proof}

It follows at once from~\eqref{eq:cassini} that $\gcd(f_n,f_{n+1})=1$: every common divisor of $f_n$ and $f_{n+1}$ divides both terms on the left-hand side and hence divides $(-1)^n$. The quotients $f_n/f_{n+1}$ are therefore irreducible fractions, and an index shift in~\eqref{eq:cassini} gives $f_nf_{n+2}-f_{n+1}^2=(-1)^{n+1}$. Hence two consecutive quotients have determinant $\pm1$ and are Farey neighbours in this sense.

\section{The metallic circle}\label{sec:metallic}

By Cassini's identity~\eqref{eq:cassini}, $\det\left( \begin{smallmatrix} f_{n-1} & f_{n} \\ f_{n} &f_{n+1} \end{smallmatrix} \right)=(-1)^{n}$, so by Proposition~\ref{prop:ford} the consecutive Ford circles $K[f_{n-1},f_{n}]$ and $K[f_{n},f_{n+1}]$ are tangent.\footnote{If some $f_n$ is negative, the signs of numerator and denominator are changed, which affects neither the criterion of Proposition~\ref{prop:ford} nor formula~\eqref{eq:tangencyK}; when $f_n=0$ the circle $K[1,0]$ is used.} By Proposition~\ref{prop:tangency} and formulas~\eqref{eq:add1} and~\eqref{eq:add2}, their tangency point is
\begin{equation*}
    \left( \frac{f_{n-1}f_n + f_n f_{n+1}}{f_n^{2} + f_{n+1}^{2}}, \frac{1}{f_n^{2} + f_{n+1}^{2}} \right)=
\left(\frac{f_{2n}}{f_{2n+1}},\;\frac{1}{f_{2n+1}} \right)=\mathcal{P}_{n},
\end{equation*}
with $n \in \mathbb{Z}$. Note that $f_{2n+1}>0$, since it is a sum of squares; hence the tangency points do lie in the upper half-plane.

\begin{theorem}[Metallic circle]\label{thm:metalliccircle}
For each integer $k\geqslant 1$, the points
\[\mathcal{P}_{n}=\left(\frac{f_{2n}}{f_{2n+1}},\;\frac{1}{f_{2n+1}} \right), \qquad n\in\mathbb{Z},\]
lie on the circle
\begin{equation}\label{eq:MetalCircle}
    C_{k}:\qquad \left( x+ \frac{k}{2}\right)^{2} + y^{2}=\frac{k^{2}+4}{4},
\end{equation}
with centre $\left(-\frac{k}{2},0\right)$ and radius $r_{k}=\frac{1}{2}\sqrt{k^{2}+4}$. Moreover, the circle meets the $X$-axis at the points $(\sigma_{k}^{-1},0)$ and $(-\sigma_{k},0)$.
\end{theorem}

\begin{proof}
Let $x=\frac{f_{2n}}{f_{2n+1}}$ and $y=\frac{1}{f_{2n+1}}$. Take identity~\eqref{eq:inv} of Lemma~\ref{lem:identities} with $n$ replaced by $2n$:
\begin{equation*}
        f_{2n+1}^2-k\,f_{2n}f_{2n+1}-f_{2n}^2=(-1)^{2n}=1.
\end{equation*}
Dividing by $f_{2n+1}^{2}$,
\[
1 - k\,\frac{f_{2n}}{f_{2n+1}}
    - \frac{f_{2n}^{2}}{f_{2n+1}^{2}}
    = \frac{1}{f_{2n+1}^{2}},
\]
which in terms of $x$ and $y$ reads
\[1 - kx - x^{2} = y^{2}\;\iff\;x^{2} + kx + y^{2} = 1.\]
Completing the square gives
\begin{align*}
    x^{2}+kx+\frac{k^{2}}{4}-\frac{k^{2}}{4}+y^{2}&=1\\
    \left( x^{2}+kx+\frac{k^{2}}{4}\right)+y^{2}&=1+\frac{k^{2}}{4} \\
\left( x+ \frac{k}{2}\right)^{2} + y^{2}&=\frac{k^{2}+4}{4},
\end{align*}
which is precisely the equation of $C_{k}$. Since identity~\eqref{eq:inv} holds for all $n\in \mathbb{Z}$, we conclude that $\mathcal{P}_{n} \in C_{k}$.

It remains to show that the intersections with the $X$-axis are $(\sigma_{k}^{-1},0)$ and $(-\sigma_{k},0)$. To this end, set $y=0$ in~\eqref{eq:MetalCircle}:
\begin{align*}
    \left( x + \frac{k}{2} \right)^{2} = \frac{k^{2} + 4}{4} &\Longleftrightarrow x^{2} + kx = 1 \\
    &\Longleftrightarrow x^{2} + kx - 1=0,
\end{align*}
whose roots are
\[
 x_{+} = \frac{-k + \sqrt{k^{2}+4}}{2} \quad\text{and}\quad x_{-} = \frac{-k - \sqrt{k^{2}+4}}{2}.
\]
Since $x_{+}=-\sigma_k'=\sigma_k^{-1}$ and $x_{-}=-\sigma_k$, these are exactly the required points.
\end{proof}

\begin{figure}[t]
\centering
\begin{tikzpicture}[scale=2.8]
  \begin{scope}
  \clip (-2.8,0) rectangle (1.52,1.6);
  \draw[->,gray] (-2.75,0) -- (1.35,0);
  \draw[violet!80!black,thick] (-2.75,1) -- (1.35,1);
  \begin{scope}[gray!40]
    \foreach \p in {-2,-1,0,1,2}{\draw (\p,0.5) circle (0.5);}
    \foreach \x in {-2.5,-1.5,-0.5,0.5,1.5}{\draw (\x,0.125) circle (0.125);}
    \foreach \x in {-1.6667,-1.3333,-0.6667,-0.3333,0.3333,0.6667,1.3333,1.6667}{\draw (\x,0.0556) circle (0.0556);}
  \end{scope}
  \begin{scope}[blue!70!black,thick]
    \draw (-1.5,0.125) circle (0.125);
    \draw (-1,0.5) circle (0.5);
    \draw (1,0.5) circle (0.5);
    \draw (0.6667,0.0556) circle (0.0556);
  \end{scope}
  \begin{scope}[red!70!black,thick]
    \draw (-2.5,0.125) circle (0.125);
    \draw (0.4,0.02) circle (0.02);
  \end{scope}
  \begin{scope}[violet!80!black,thick]
    \draw (-2,0.5) circle (0.5);
    \draw (0,0.5) circle (0.5);
    \draw (0.5,0.125) circle (0.125);
  \end{scope}
  \draw[blue!70!black,very thick] (-0.5,0) circle (1.118034);
  \draw[red!70!black,very thick] (-1,0) circle (1.414214);
  \end{scope}
  \node[below right,gray] at (1.35,0) {\small $X$};
  \node[right,violet!80!black] at (1.35,1) {\small $y=1$};
  \foreach \pt in {(-1.6,0.2),(-1.5,0.5),(-1,1),(0.5,0.5),(0.6,0.2),(0.615385,0.076923)}{\fill[blue!70!black] \pt circle (0.028);}
  \foreach \pt in {(-2.4,0.2),(-2,1),(0.4,0.2),(0.413793,0.034483)}{\fill[red!70!black] \pt circle (0.024);}
  \fill[violet!80!black] (0,1) circle (0.032);
  \fill (0.618034,0) circle (0.022);
  \node[below right,inner sep=1pt] at (0.618034,-0.02) {\scriptsize $\varphi^{-1}$};
  \fill (-1.618034,0) circle (0.022);
  \node[below] at (-1.618034,-0.02) {\scriptsize $-\varphi$};
  \fill (0.414214,0) circle (0.022);
  \node[below left,red!60!black,inner sep=1pt] at (0.414214,-0.02) {\scriptsize $\delta_{\text{Ag}}^{-1}$};
  \fill (-2.414214,0) circle (0.022);
  \node[below,red!60!black] at (-2.414214,-0.02) {\scriptsize $-\delta_{\text{Ag}}$};
  \node[blue!70!black] at (-0.5,1.22) {\small $C_1$};
  \node[red!70!black] at (-1.0,1.49) {\small $C_2$};
  \node[violet!80!black] at (0.13,0.88) {\scriptsize $\mathcal{P}_0$};
\end{tikzpicture}
\caption{The metallic circles $C_1$ (golden, blue) and $C_2$ (silver, red), with the
chains of Ford circles $K[f_{n-1},f_n]$ that generate them; in violet, the elements
common to both chains and the point $\mathcal{P}_0=(0,1)$. Each marked point is the
tangency point of two consecutive circles of its chain. $C_1$ meets the axis at $\varphi^{-1}$
and $-\varphi$; $C_2$, at $\delta_{\text{Ag}}^{-1}$ and $-\delta_{\text{Ag}}$. The case $C_1$
is Kocik's.}
\label{fig:metallic}
\end{figure}
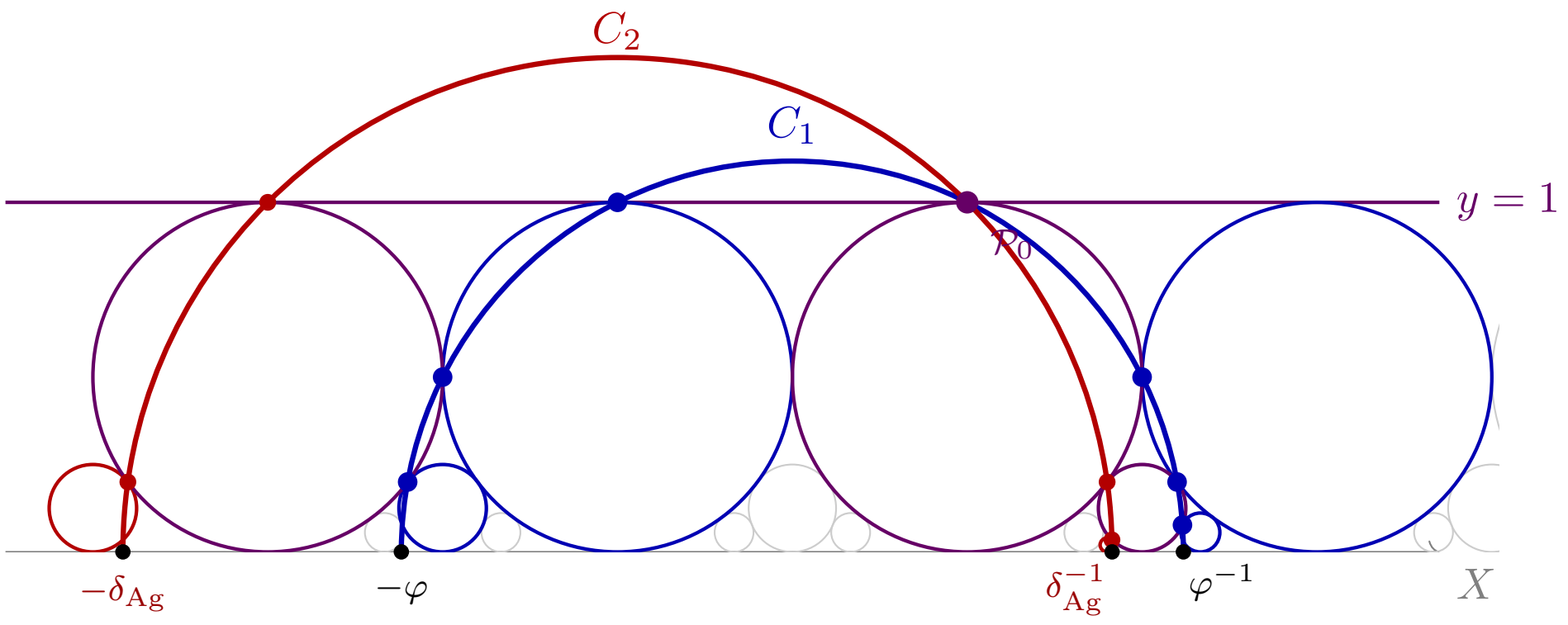

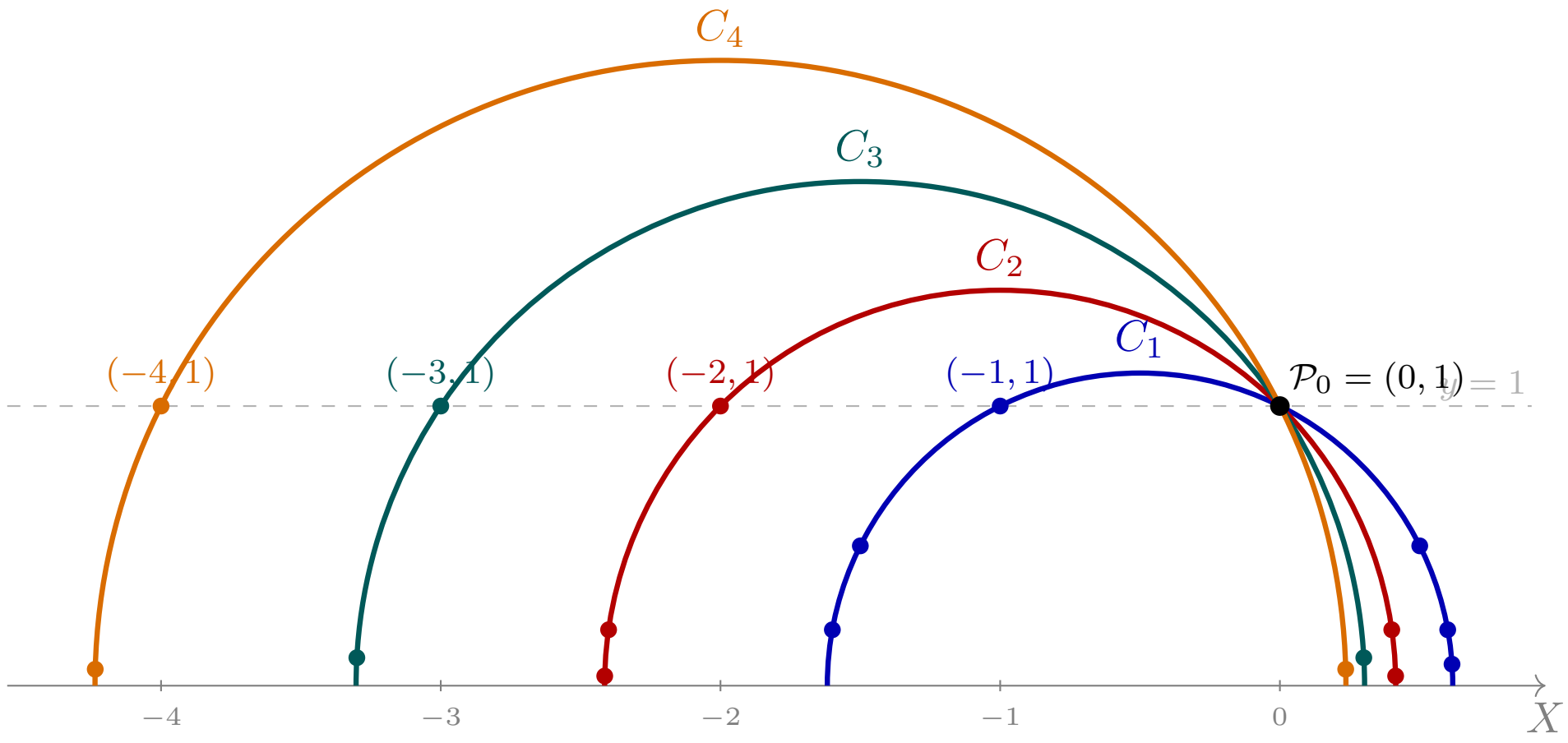
\begin{figure}[t]
\centering
\begin{tikzpicture}[scale=2.25]
  \draw[->,gray] (-4.55,0) -- (0.95,0) node[below,gray] {\small $X$};
  \draw[gray!60,dashed] (-4.55,1) -- (0.9,1) node[above left,inner sep=1pt] {\scriptsize $y=1$};
  \draw[blue!70!black,very thick]  (0.618034,0) arc (0:180:1.118034);
  \draw[red!70!black,very thick]   (0.414214,0) arc (0:180:1.414214);
  \draw[teal!70!black,very thick]  (0.302776,0) arc (0:180:1.802776);
  \draw[orange!85!black,very thick](0.236068,0) arc (0:180:2.236068);
  \foreach \pt in {(-1.6,0.2),(-1.5,0.5),(-1,1),(0,1),(0.5,0.5),(0.6,0.2),(0.6154,0.0769)}{\fill[blue!70!black] \pt circle (0.03);}
  \foreach \pt in {(-2.4138,0.0345),(-2.4,0.2),(-2,1),(0,1),(0.4,0.2),(0.4138,0.0345)}{\fill[red!70!black] \pt circle (0.03);}
  \foreach \pt in {(-3.3,0.1),(-3,1),(0,1),(0.3,0.1)}{\fill[teal!70!black] \pt circle (0.03);}
  \foreach \pt in {(-4.2353,0.0588),(-4,1),(0,1),(0.2353,0.0588)}{\fill[orange!85!black] \pt circle (0.03);}
  \fill[black] (0,1) circle (0.035);
  \node[above right,inner sep=2pt] at (0,1) {\scriptsize $\mathcal P_0=(0,1)$};
  \node[above,blue!70!black,inner sep=3pt] at (-1,1) {\scriptsize $(-1,1)$};
  \node[above,red!70!black,inner sep=3pt] at (-2,1) {\scriptsize $(-2,1)$};
  \node[above,teal!70!black,inner sep=3pt] at (-3,1) {\scriptsize $(-3,1)$};
  \node[above,orange!85!black,inner sep=3pt] at (-4,1) {\scriptsize $(-4,1)$};
  \node[blue!70!black]  at (-0.5,1.24) {\small $C_1$};
  \node[red!70!black]   at (-1.0,1.53) {\small $C_2$};
  \node[teal!70!black]  at (-1.5,1.92) {\small $C_3$};
  \node[orange!85!black]at (-2.0,2.35) {\small $C_4$};
  \foreach \x in {-4,-3,-2,-1,0}{\draw[gray] (\x,0.02) -- (\x,-0.02) node[below,gray] {\tiny $\x$};}
\end{tikzpicture}
\caption{The first four metallic circles $C_1$ (golden), $C_2$ (silver),
$C_3$ and $C_4$, with the tangency points $\mathcal P_n$ visible at this scale. All of them
pass through $\mathcal P_0=(0,1)$, and the point $\mathcal P_{-1}$ of $C_k$ is $(-k,1)$. As
$k$ grows, the right foot $\sigma_k^{-1}$ approaches $0$ and the left foot $-\sigma_k$
approaches $-k$, while the points $\mathcal P_n$ crowd ever faster towards the feet.}
\label{fig:family}
\end{figure}

Note that for $k=1$ we recover Kocik's result: the circle with centre $\left( -\frac{1}{2},0 \right)$ and radius $\frac{\sqrt{5}}{2}$, meeting the $X$-axis at $(\varphi^{-1},0)$ and $(-\varphi,0)$. For $k=2$, the circle associated with the Pell sequence has centre $(-1,0)$ and radius $\sqrt{2}$, and meets the $X$-axis at $(-1+\sqrt{2},0)$ and $(-1-\sqrt{2},0)$, that is, at $\delta_{\text{Ag}}^{-1}$ and $-\delta_{\text{Ag}}$, respectively (see \cref{fig:metallic}). \Cref{fig:family} shows the first four circles of the family.

We call the two points where $C_{k}$ meets the $X$-axis the \textit{feet} of $C_{k}$.

\begin{corollary}[The feet are the metallic numbers]\label{cor:feet}
The sequence of tangency points runs along $C_{k}$ from one foot to the other: the abscissa $x_n=f_{2n}/f_{2n+1}$ of $\mathcal{P}_n$ is strictly increasing in $n \in \mathbb{Z}$, and
\begin{equation}
\lim_{n\to+\infty}\mathcal{P}_n=\left(\sigma_k^{-1},0\right),  \qquad
  \lim_{n\to-\infty}\mathcal{P}_n=\left(-\sigma_k,0\right).
\end{equation}
\end{corollary}

\begin{proof}
Since $f_{2n}/f_{2n+1}\to \sigma_{k}^{-1}$ and $1/f_{2n+1} \to 0$ as $n \to +\infty$, we have $\mathcal{P}_{n}\to ( \sigma_{k}^{-1},0)$. For $n\to -\infty$, using $f_{-m}=(-1)^{m+1}f_{m}$ we get $f_{2n}/f_{2n+1} = - f_{2|n|}/f_{2|n|-1}\to -\sigma_{k}$ and $1/f_{2|n|-1} \to 0^{+}$, that is, $\mathcal{P}_{n} \to (-\sigma_{k},0)$.

For monotonicity, d'Ocagne's identity~\eqref{eq:docagne} with $a=2n+2$ and $b=2n$ gives $f_{2n+2}f_{2n+1}-f_{2n+3}f_{2n}=f_2=k$. Hence
\begin{equation*}
x_{n+1}-x_n =\frac{f_{2n+2}f_{2n+1}-f_{2n+3}f_{2n}}{f_{2n+1}f_{2n+3}}=\frac{k}{f_{2n+1}f_{2n+3}}>0 ,
\end{equation*}
since, by~\eqref{eq:add2}, both factors in the denominator are positive for all $n\in \mathbb{Z}$. Thus $\{x_{n}\}$ increases strictly from $-\sigma_k$ to $\sigma_k^{-1}$ and, since the points lie in the upper half-plane, the $\mathcal{P}_{n}$ traverse the upper semicircle of $C_{k}$ from one foot to the other.
\end{proof}

\begin{proposition}[The translated and reflected chain]\label{prop:otherchains}
Let $j\in\mathbb{Z}$. The translation $(x,y)\mapsto (x+j,y)$ and the reflection
$(x,y)\mapsto (-x,y)$ preserve the family of Ford circles and hence map
tangency points to tangency points. In particular, the images of the points
$\mathcal{P}_{n}$ under the translation with $j=k$
\begin{align*}
   &\qquad \left(\frac{f_{2n}}{f_{2n+1}}+k,\ \frac{1}{f_{2n+1}}\right)
       =\left(\frac{f_{2n+2}}{f_{2n+1}},\ \frac{1}{f_{2n+1}}\right),\\[2pt]
  &\text{and under the reflection:}\qquad \left(-\frac{f_{2n}}{f_{2n+1}},\ \frac{1}{f_{2n+1}}\right),
\end{align*}
are tangency points of Ford chains whose coordinates are quotients of terms
of the $k$-Fibonacci sequence. Both families coincide as sets (the reflected point
with index $n$ is the translated point with index $-n-1$) and they lie on the circle
with centre $\left(\frac{k}{2},0\right)$ and radius $r_{k}$, which meets the axis at
$(\sigma_k,0)$ and $(-\sigma_k^{-1},0)$.
\end{proposition}

\begin{proof}
The first equality is the recurrence $f_{2n+2}=kf_{2n+1}+f_{2n}$. Since
$f_{-2n}=-f_{2n}$ and $f_{-2n-1}=f_{2n+1}$, the translated point with index $-n-1$ is
\[
  \left(\frac{f_{-2n}}{f_{-2n-1}},\ \frac{1}{f_{-2n-1}}\right)
  =\left(-\frac{f_{2n}}{f_{2n+1}},\ \frac{1}{f_{2n+1}}\right),
\]
that is, the reflected point with index $n$. Finally, if $(x,y)\in C_{k}$ then
$x^{2}+kx+y^{2}=1$, so $(-x,y)$ satisfies $x^{2}-kx+y^{2}=1$, which is the
equation of the circle with centre $\left(\frac k2,0\right)$ and radius $r_{k}$. The
intersections with the axis are obtained as in Theorem~\ref{thm:metalliccircle}.
\end{proof}

\begin{remark}
The coincidence of both families reflects the fact that $C_{k}$ is symmetric about the
line $x=-\frac{k}{2}$: the symmetry $(x,y)\mapsto(-k-x,y)$ swaps
$\mathcal{P}_{n}$ and $\mathcal{P}_{-n-1}$. For $k=1$ we recover exactly Kocik's
chain (B) \citep{kocik2020fibonacci}. As for Kocik's chain (C), which for $k=1$ is obtained by translating by $2$
and has coordinates $F_{2n+3}/F_{2n+1}$, the analogous translation
$(x,y)\mapsto (x+2k,y)$ produces $\left(f_{2n+2}+kf_{2n+1}\right)/f_{2n+1}$, and
$f_{2n+2}+kf_{2n+1}$ is not in general a term of the sequence (it is only when $k=1$).
Whether a metallic analogue of (C) and (D) exists is left as an open question.
\end{remark}

\begin{figure}[t]
\centering
\begin{tikzpicture}[scale=1.9]
  \begin{scope}
  \clip (-2.85,0) rectangle (2.85,1.3);
  \draw[->,gray] (-2.8,0) -- (2.8,0);
  \draw[gray!60,dashed] (-2.8,1) -- (2.8,1);
  \begin{scope}[gray!35]
    \foreach \p in {-3,-2,-1,0,1,2,3}{\draw (\p,0.5) circle (0.5);}
    \foreach \x in {-2.5,-1.5,-0.5,0.5,1.5,2.5}{\draw (\x,0.125) circle (0.125);}
  \end{scope}
  \draw[orange!85!black,very thick] (-1.5,0) circle (1.118034);
  \draw[blue!70!black,very thick]   (-0.5,0) circle (1.118034);
  \draw[red!70!black,very thick]    ( 0.5,0) circle (1.118034);
  \draw[teal!70!black,very thick]   ( 1.5,0) circle (1.118034);
  \end{scope}
  \foreach \j/\c in {-1/orange!85!black,0/blue!70!black,1/red!70!black,2/teal!70!black}{
    \foreach \x/\y in {-1.6/0.2,-1.5/0.5,-1/1,0/1,0.5/0.5,0.6/0.2}{
      \fill[\c] ({\x+\j},\y) circle (0.035);}}
  \node[orange!85!black] at (-1.5,1.22) {\small (D)};
  \node[blue!70!black]   at (-0.5,1.22) {\small (A)};
  \node[red!70!black]    at ( 0.5,1.22) {\small (B)};
  \node[teal!70!black]   at ( 1.5,1.22) {\small (C)};
  \foreach \x in {-2,-1,0,1,2}{\draw[gray] (\x,0.02) -- (\x,-0.02) node[below,gray] {\tiny $\x$};}
\end{tikzpicture}
\caption{Kocik's four concyclic sequences for $k=1$: (A) consists of the points
$\mathcal P_n$, on $C_1$; (B), (C) and (D) are obtained by translating (A) by $1$, $2$ and
$-1$, respectively. The four circles have radius $\frac{\sqrt5}{2}$ and centres
$-\frac12,\frac12,\frac32,-\frac32$. For general $k$, the translation by $k$ and the
reflection produce the same family (Proposition~\ref{prop:otherchains}).}
\label{fig:kocikABCD}
\end{figure}
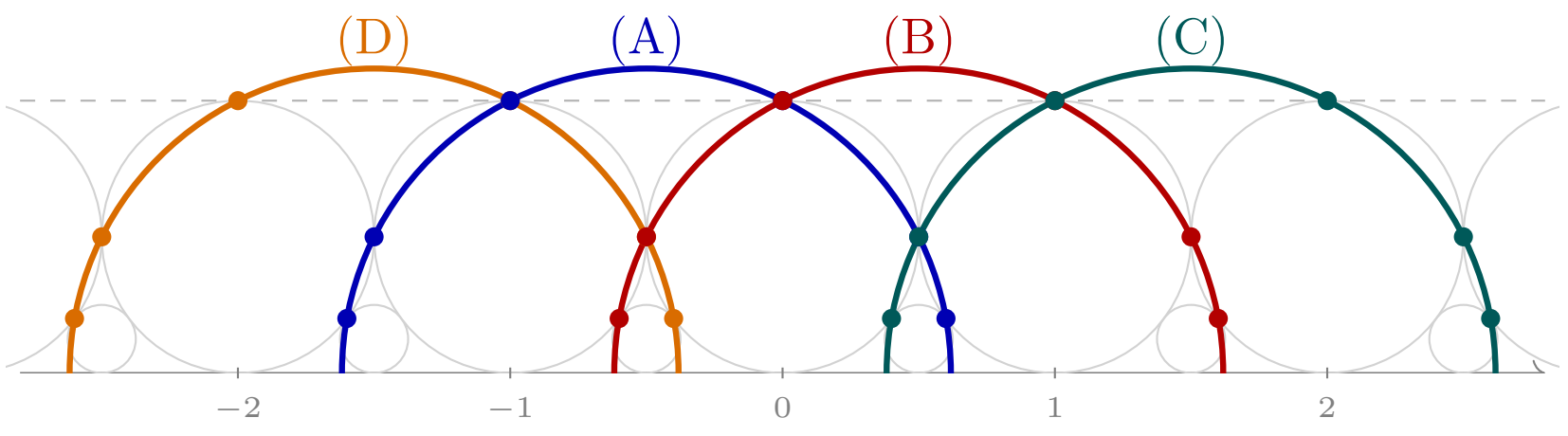

\Cref{fig:kocikABCD} shows Kocik's four chains for $k=1$.

\subsection{A hyperbolic reading}\label{sec:hyperbolic}

The proof of Theorem~\ref{thm:metalliccircle} is purely arithmetic, but it is worth pointing out where the phenomenon comes from. Identify $(x,y)$ with $z=x+iy$ in the upper half-plane $\mathbb{H}$, endowed with the hyperbolic metric, whose distance is given by
\[
\cosh d_{\mathbb H}(z,w)=1+\frac{|z-w|^2}{2\,\operatorname{Im}z\,\operatorname{Im}w}.
\]
If $g=\left(\begin{smallmatrix} p & r \\ q & s \end{smallmatrix}\right)$ has $\det g=ps-qr=1$, then
\begin{equation}\label{eq:action}
  g\cdot i=\frac{pi+r}{qi+s}=\frac{pq+rs}{q^{2}+s^{2}}+\frac{i}{q^{2}+s^{2}},
\end{equation}
which is exactly the tangency formula~\eqref{eq:tangencyK}. Recall also that $\mathrm{SL}_2(\mathbb Z)$, acting by Möbius transformations, and the reflection $z\mapsto-\bar z$ are isometries of $\mathbb H$ that permute the Ford circles, including $K[1,0]$ \citep{ford1938fractions}. Hence every element of $\mathrm{GL}_2(\mathbb Z)$ acts in this way: those of determinant $1$ by $z\mapsto g\cdot z$, and those of determinant $-1$ by $z\mapsto g\cdot\bar z$.

Now let
\[
N=\begin{pmatrix} 0 & 1 \\ 1 & k\end{pmatrix}=JMJ,\qquad
J=\begin{pmatrix} 0 & 1 \\ 1 & 0\end{pmatrix},
\]
so that $N^{n}=JM^{n}J=\left(\begin{smallmatrix} f_{n-1} & f_{n} \\ f_{n} &
f_{n+1}\end{smallmatrix}\right)$ and $\det N^{n}=(-1)^{n}$. Since $\det N=-1$, the matrix
$N$ acts on $\mathbb{H}$ by the antiholomorphic isometry
\[
  \nu(z)=\frac{1}{\bar z+k},
\]
whose square is the Möbius transformation associated with
$N^{2}\in\mathrm{SL}_{2}(\mathbb{Z})$. The map $\nu$ is the composition of the
translation $z\mapsto z+k$, the reflection $z\mapsto -\bar z$ and the inversion
$z\mapsto -1/z$, all of which preserve the family of Ford circles; moreover, it maps the
fraction $f_{n-1}/f_{n}$ to $f_{n}/f_{n+1}$ and hence the tangency point
$\mathcal{P}_{n}$ to $\mathcal{P}_{n+1}$. Since $\mathcal{P}_{0}=i$, we obtain
\begin{equation}\label{eq:orbit}
  \mathcal{P}_{n}=\nu^{n}(i),\qquad \mathcal{P}_{2n}=N^{2n}\cdot i,\qquad n\in\mathbb{Z},
\end{equation}
where the second equality can also be read off directly from~\eqref{eq:action},
\eqref{eq:add1} and~\eqref{eq:add2}. In other words, the tangency points form the
orbit of $i$ under $\nu$, and those with even index form the orbit of $i$ under the
element $N^{2}$. Since
$\operatorname{tr}N^{2}=k^{2}+2>2$, this element is hyperbolic, and its fixed points are those of $\nu$ on the real axis, namely the roots of $z^{2}+kz-1=0$: precisely $\sigma_k^{-1}$ and $-\sigma_k$. Its axis ---the geodesic of $\mathbb{H}$ joining both fixed points--- is the semicircle with centre $\left(-\frac k2,0\right)$ and radius $\frac{\sqrt{k^{2}+4}}{2}$, that is, the upper half of $C_{k}$. Since $\nu$ is an isometry fixing both endpoints, it preserves this axis, which contains $i$. In the usual terminology, $\nu$ is a glide reflection: it composes the reflection in the geodesic $C_k\cap\mathbb H$ with a translation along it. The following lemma computes the length of that translation.

\begin{lemma}[Distance between consecutive tangency points]\label{lem:dist}
Let $p/q$ and $r/s$ be fractions with $ps-qr=\pm1$, let $a\geqslant1$ be an integer and let $(u,v)=a(r,s)+(p,q)$. If $T$ is the tangency point of $K[p,q]$ and $K[r,s]$, and $T'$ that of $K[r,s]$ and $K[u,v]$, then
\[
d_{\mathbb H}(T,T')=2\operatorname{arcsinh}\frac a2.
\]
\end{lemma}

\begin{proof}
Let $h=\left(\begin{smallmatrix}r&p\\s&q\end{smallmatrix}\right)\in\mathrm{GL}_2(\mathbb Z)$, which acts as an isometry of $\mathbb H$ permuting the Ford circles. Since $h(\infty)=r/s$, $h(0)=p/q$ and $h(a)=u/v$, it maps $K[1,0]$, $K[0,1]$ and $K[a,1]$ to $K[r,s]$, $K[p,q]$ and $K[u,v]$, respectively, and therefore maps tangency points to tangency points. The line $K[1,0]=\{y=1\}$ touches $K[0,1]$ at $i$ and $K[a,1]$ at $a+i$. Hence
\[
\cosh d_{\mathbb H}(T,T')=\cosh d_{\mathbb H}(i,\,a+i)=1+\frac{a^2}{2},
\]
and since $\cosh(2u)=1+2\sinh^2u$, we obtain $d_{\mathbb H}(T,T')=2\operatorname{arcsinh}(a/2)$.
\end{proof}

\begin{corollary}[Equidistance]\label{cor:equidistance}
The points $\mathcal P_n$ are equally spaced along the geodesic $C_k\cap\mathbb H$:
\[
d_{\mathbb H}(\mathcal P_m,\mathcal P_n)=2\,|m-n|\,\ln\sigma_k,\qquad m,n\in\mathbb Z.
\]
In particular, $\nu$ is a glide reflection with translation length $2\ln\sigma_k$.
\end{corollary}

\begin{proof}
By the recurrence, $(f_n,f_{n+1})=k\,(f_{n-1},f_n)+(f_{n-2},f_{n-1})$, and by~\eqref{eq:cassini} the fractions $f_{n-2}/f_{n-1}$ and $f_{n-1}/f_n$ satisfy the hypothesis of Lemma~\ref{lem:dist} with $a=k$. Hence $d_{\mathbb H}(\mathcal P_{n-1},\mathcal P_n)=2\operatorname{arcsinh}(k/2)$, and
\[
\operatorname{arcsinh}\frac k2=\ln\left(\frac k2+\sqrt{\frac{k^2}{4}+1}\right)=\ln\sigma_k.
\]
Since the $\mathcal P_n$ lie on a geodesic and advance monotonically (Corollary~\ref{cor:feet}), the distances between consecutive points add up.
\end{proof}

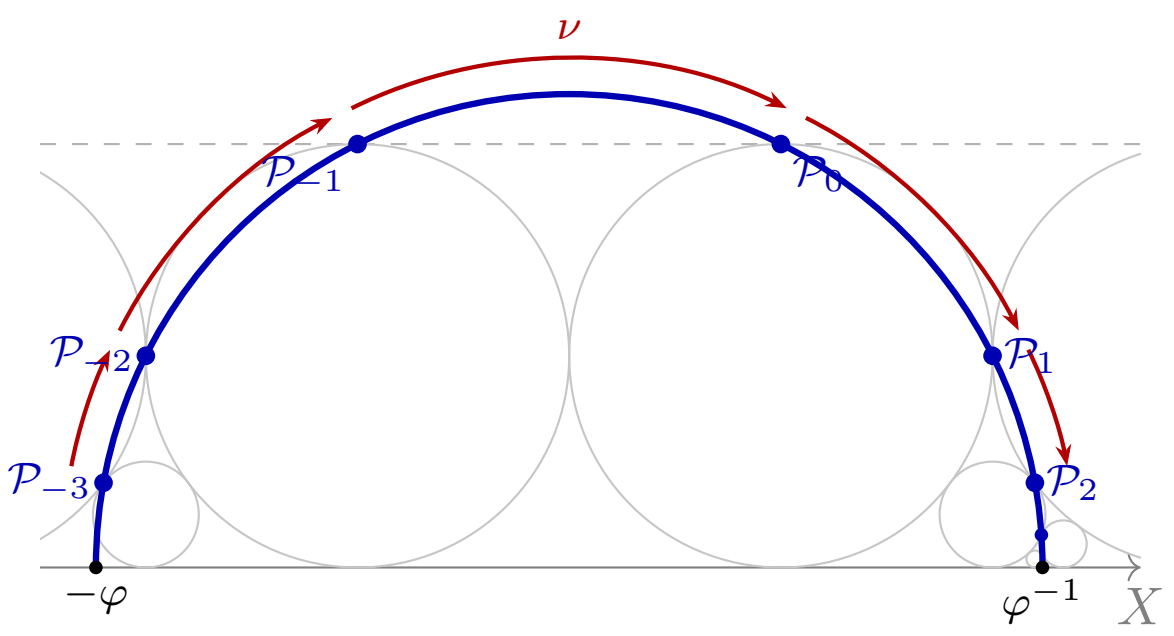
\begin{figure}[htbp]
\centering
\begin{tikzpicture}[scale=2.8]
  \draw[->,gray] (-1.75,0) -- (0.85,0) node[below,gray] {\small $X$};
  \draw[gray!60,dashed] (-1.75,1) -- (0.85,1);
  \begin{scope}
  \clip (-1.75,0) rectangle (0.85,1.3);
  \begin{scope}[gray!45]
    \draw (-1.5,0.125) circle (0.125);
    \draw (-2,0.5) circle (0.5);
    \draw (-1,0.5) circle (0.5);
    \draw (0,0.5) circle (0.5);
    \draw (1,0.5) circle (0.5);
    \draw (0.5,0.125) circle (0.125);
    \draw (0.6667,0.0556) circle (0.0556);
    \draw (0.6,0.02) circle (0.02);
  \end{scope}
  \end{scope}
  \draw[blue!70!black,very thick] (0.618034,0) arc (0:180:1.118034);
  \foreach \a/\b in {169.7/153.43,153.43/116.57,116.57/63.43,63.43/26.57,26.57/10.3}{
    \draw[-{Stealth[length=4pt]},red!70!black,thick,shorten <=2pt,shorten >=2pt]
      ([shift={(\a:1.20)}]-0.5,0) arc (\a:\b:1.20);}
  \node[red!70!black,above] at ([shift={(90:1.20)}]-0.5,0) {\scriptsize $\nu$};
  \foreach \x/\y/\n/\pos in {-1.6/0.2/-3/left,-1.5/0.5/-2/left,-1/1/-1/below left,0/1/0/below right,0.5/0.5/1/right,0.6/0.2/2/right}{
    \fill[blue!70!black] (\x,\y) circle (0.022);
    \node[\pos,inner sep=2pt,blue!70!black] at (\x,\y) {\scriptsize $\mathcal P_{\n}$};}
  \fill[blue!70!black] (0.6154,0.0769) circle (0.016);
  \fill (0.618034,0) circle (0.016) node[below,inner sep=2pt] {\scriptsize $\varphi^{-1}$};
  \fill (-1.618034,0) circle (0.016) node[below,inner sep=2pt] {\scriptsize $-\varphi$};
\end{tikzpicture}
\caption{The geodesic $C_1\cap\mathbb H$ and Kocik's points $\mathcal P_n$,
with the Fibonacci Ford chain in the background. The glide reflection $\nu$ maps each
$\mathcal P_n$ to $\mathcal P_{n+1}$ (arrows). Although in the Euclidean picture the points
crowd towards the feet, in the hyperbolic metric all steps have the same length,
$2\ln\varphi\approx0.962$ (Corollary~\ref{cor:equidistance}).}
\label{fig:hyperbolic}
\end{figure}

For instance, Kocik's points are separated by $2\ln\varphi\approx0.962$, and those of the Pell chain by $2\ln(1+\sqrt2)\approx1.763$ (see \cref{fig:hyperbolic}).

From this point of view, Theorem~\ref{thm:metalliccircle} says that the orbit of a point on the axis of a hyperbolic element stays on that axis, and Corollary~\ref{cor:feet} is nothing but the north--south dynamics of that element: the points advance monotonically from the repelling fixed point $-\sigma_k$ towards the attracting one $\sigma_k^{-1}$. Corollary~\ref{cor:equidistance} makes the speed precise: a step of constant length $2\ln\sigma_k$. This explains \textit{why} concyclicity is not accidental.

\section{Concluding remarks}\label{sec:conclusions}

Kocik's circle is not a peculiarity of the golden ratio: every metallic ratio $\sigma_k$ has its own, and its origin is arithmetic rather than geometric. Concyclicity is identity~\eqref{eq:inv} read in coordinates. Hyperbolic geometry adds the scale: the tangency points advance along the geodesic $C_k$ with constant step $2\ln\sigma_k$ (Corollary~\ref{cor:equidistance}).



\begin{thebibliography}{7}
\providecommand{\natexlab}[1]{#1}
\providecommand{\url}[1]{\texttt{#1}}
\expandafter\ifx\csname urlstyle\endcsname\relax
  \providecommand{\doi}[1]{doi: #1}\else
  \providecommand{\doi}{doi: \begingroup \urlstyle{rm}\Url}\fi

\bibitem[Apostol(1990)]{apostol1997modular}
Tom~M. Apostol.
\newblock \emph{Modular Functions and {Dirichlet} Series in Number Theory}.
\newblock Number~41 in Graduate Texts in Mathematics. Springer-Verlag, New
  York, second edition, 1990.
\newblock \doi{10.1007/978-1-4612-0999-7}.

\bibitem[{de Spinadel}(1999)]{spinadel1999metallic}
Vera~W. {de Spinadel}.
\newblock The family of metallic means.
\newblock \emph{Visual Mathematics}, 1\penalty0 (3), 1999.
\newblock URL \url{http://geodesic.mathdoc.fr/item/VM_1999_1_3_a3/}.

\bibitem[Falc{\'o}n and Plaza(2007)]{falcon2007kfibonacci}
Sergio Falc{\'o}n and {\'A}ngel Plaza.
\newblock On the {Fibonacci} $k$-numbers.
\newblock \emph{Chaos, Solitons \& Fractals}, 32\penalty0 (5):\penalty0
  1615--1624, 2007.
\newblock \doi{10.1016/j.chaos.2006.09.022}.

\bibitem[Ford(1938)]{ford1938fractions}
Lester~R. Ford.
\newblock Fractions.
\newblock \emph{The American Mathematical Monthly}, 45\penalty0 (9):\penalty0
  586--601, 1938.
\newblock \doi{10.2307/2302799}.

\bibitem[Kocik(2020)]{kocik2020fibonacci}
Jerzy Kocik.
\newblock {Fibonacci} numbers and {Ford} circles, 2020.
\newblock URL \url{https://arxiv.org/abs/2003.00852}.
\newblock arXiv:2003.00852.

\bibitem[Sporn(2022)]{Sporn2022}
Howard Sporn.
\newblock {Fibonacci} fraction circles.
\newblock \emph{The Mathematical Gazette}, 106\penalty0 (565):\penalty0 1--8,
  2022.
\newblock \doi{10.1017/mag.2022.2}.

\bibitem[Vajda(1989)]{vajda1989fibonacci}
Steven Vajda.
\newblock \emph{{Fibonacci} and {Lucas} Numbers, and the Golden Section: Theory
  and Applications}.
\newblock Ellis Horwood, Chichester, 1989.
\newblock Reprinted by Dover, Mineola, NY, 2008.

\end{thebibliography}
\end{document}